\documentclass{aims-ppn}
\usepackage{amsmath}
\usepackage{paralist}
\usepackage{graphics}
\usepackage{epsfig}
\usepackage{graphicx}
\usepackage{epstopdf}
\usepackage[colorlinks=true]{hyperref}
\hypersetup{urlcolor=blue, citecolor=red}

\def\currentvolume{X}

\def\ppages{X--XX}

\theoremstyle{definition}

\usepackage{setspace,xspace,graphicx}
\usepackage[colorlinks=true]{hyperref}
\hypersetup{urlcolor=blue, citecolor=red, linkcolor= blue}
\usepackage{cancel}
\usepackage{color}

\newtheorem{Thm}{Theorem}
\newtheorem{Cor}[Thm]{Corollary}

\newtheorem{Prop}[Thm]{Proposition}
\newtheorem{Rem}[Thm]{Remark}

\newcommand{\R}{{\mathbb R}}
\renewcommand{\S}{{\mathbb S}}

\renewcommand{\L}{\mathrm L}
\renewcommand{\H}{\mathrm H}
\newcommand{\be}[1]{\begin{equation}\label{#1}}
\newcommand{\ee}{\end{equation}}
\renewcommand{\(}{\left(}
\renewcommand{\)}{\right)}

\newcommand{\irdeux}[1]{\int_{\R^2}{#1}\,dx}
\newcommand{\irdmu}[1]{\int_{\R^2}{#1}\,d\mu}

\newcommand{\M}{\mathcal M}
\newcommand{\iM}[1]{\int_{\M}{#1}\,d\kern1pt v_g}
\newcommand{\nrM}[2]{\|{#1}\|_{\L^{#2}(\M)}}
\newcommand{\Lap}{\Delta}
\newcommand{\Ric}{\mathrm{Ric}\,}

\newcommand{\irdnu}[1]{\int_{\R^2}{#1}\,d\nu}
\newtheorem{Ex}{Example}

\title[Onofri inequalities and rigidity results]{Onofri inequalities and rigidity results}
\author[Jean Dolbeault and Maria J.~Esteban and Gaspard Jankowiak]{}



\email{dolbeaul@ceremade.dauphine.fr}
\email{esteban@ceremade.dauphine.fr}
\email{gaspard.jankowiak@math.crns.fr}
\thanks{$^*$ Corresponding author: Jean Dolbeault}

\begin{document}
\maketitle
\thispagestyle{empty}

\centerline{\scshape Jean Dolbeault and Maria J.~Esteban$^*$}
\medskip
{\footnotesize
\centerline{Ceremade, CNRS UMR n$^{\circ}$ 7534 and Universit\'e Paris-Dauphine}
\centerline{PSL research university}
\centerline{Place de Lattre de Tassigny, 75775 Paris C\'edex~16, France}
}

\medskip

\centerline{\scshape Gaspard Jankowiak}
\medskip
{\footnotesize
\centerline{RICAM and Universit\"at Wien}
\centerline{Wolfgang Pauli Institute}
\centerline{Oskar-Morgenstern-Platz 1, 1090 Wien, Austria}
}

\bigskip

\centerline{(Communicated by Manuel del Pino)}

\begin{abstract}
This paper is devoted to the Moser-Trudinger-Onofri inequality on smooth compact connected Riemannian manifolds. We establish a rigidity result for the Euler-Lagrange equation and deduce an estimate of the optimal constant in the inequality on two-dimensional closed Riemannian manifolds. Compared to existing results, we provide a non-local criterion which is well adapted to variational methods, introduce a nonlinear flow along which the evolution of a functional related with the inequality is monotone and get an integral remainder term which allows us to discuss optimality issues. As an important application of our method, we also consider the non-compact case of the Moser-Trudinger-Onofri inequality on the two-dimensional Euclidean space, with weights. The standard weight is the one that is computed when projecting the two-dimensional sphere using the stereographic projection, but we also give more general results which are of interest, for instance, for the Keller-Segel model in chemotaxis.
\end{abstract}

\vspace*{1cm}


In this paper we assume that $(\M,g)$ is a smooth compact connected Riemannian manifold of dimension $d\ge1$, without boundary. We denote by $\Lap$ the Laplace-Beltrami operator on $\M$. For simplicity, we assume that the volume of $\M$, is chosen equal to $1$ and use the notation $d\kern1pt v_g$ for the volume element. We shall also denote by~$\Ric$ the Ricci tensor, by $\mathrm H u$ the Hessian of $u$ and by
\[
\mathrm L u:=\mathrm H u-\frac gd\,\Lap u
\]
the trace free Hessian. Let us denote by $\mathrm M[u]$ the trace free tensor
\[
\mathrm M[u]:=\nabla u\otimes\nabla u-\frac gd\,|\nabla u|^2\,.
\]
Let $g^{i,j}$ be the inverse of the metric tensor, \emph{i.e.}, $g^{i,j}\,g_{j,k}=\delta^i_k$, where $\delta^i_k$ denotes the Kronecker symbol. If $\mathrm A_{i,j}$ and $\mathrm B_{i,j}$ are two tensors, we use the notation
\[
\mathrm A:\mathrm B:=g^{i,m}\,g^{j,n}\,\mathrm A_{i,j}\,\mathrm B_{m,n}\quad\mbox{and}\quad\|\mathrm A\|^2:=\mathrm A:\mathrm A\,,
\]
where we used the Einstein summation convention. We define
\be{LambdaStar}
\lambda_\star:=\inf_{u\in\H^2(\M)\setminus\{0\}}\frac{\displaystyle\iM{\Big[\,\|\,\mathrm L u-\tfrac12\,\mathrm M[u]\,\|^2+\Ric(\nabla u,\nabla u)\Big]\,e^{-u/2}}}{\displaystyle\iM{|\nabla u|^2\,e^{-u/2}}}\,.
\ee
In this paper we will prove the two following results.
\begin{Thm}\label{Thm:Rigidity} Assume that $d=2$ and $\lambda_\star>0$. If $u$ is a smooth solution to
\be{Eqn:EL}
-\,\frac12\,\Lap u+\lambda=e^u\,,
\ee
then $u$ is a constant function if $\lambda\in(0,\lambda_\star)$.\end{Thm}
Note that the hypothesis that $\lambda_\star>0$ is in principle weaker than assuming that the Ricci curvature is everywhere positive on $\M$. See Remark~\ref{rem:positivity lambda} for more details.

Next, let us consider the Moser-Trudinger-Onofri inequality on $\M$ written as
\be{Ineq:Moser-Trudinger-Onofri}
\frac14\,\nrM{\nabla u}2^2+\lambda\,\iM u\ge\lambda\,\log\(\iM{e^u}\)\quad\forall\,u\in\H^1(\M)\,,
\ee
for some constant $\lambda>0$. Let us denote by $\lambda_1$ the first positive eigenvalue of $-\,\Lap$. We would like to draw the attention of the reader to the fact that, because of the normalization of the volume of $\mathcal M$, there is a discrepancy of a factor $4\,\pi$ with many results that are available in the literature. This factor corresponds to the surface of the usual $\S^2$ sphere, considered with the measure induced by Lebesgue's measure in $\R^3$.
\begin{Cor}\label{Cor:Rigidity} If $d=2$, then~\eqref{Ineq:Moser-Trudinger-Onofri} holds with $\lambda=\Lambda:=\min\{1,\lambda_\star\}\le \lambda_1/2$. Moreover, if $\Lambda$ is strictly smaller than $\lambda_1/2$, then the optimal constant in~\eqref{Ineq:Moser-Trudinger-Onofri} is strictly larger than~$\Lambda$.\end{Cor}
As we shall see later, in the case of the normalized sphere, $\lambda_\star=\lambda_1/2=1$ is optimal, but for $\lambda=\lambda_\star$, Eq.~\eqref{Eqn:EL} has non-constant solutions because of the conformal invariance: see~\cite{MR3052352,DEJ} for more details on the Moser-Trudinger-Onofri inequality on the sphere, and references therein. The interested reader is invited to refer to the historical papers~\cite{MR0216286,MR0301504,MR677001}, and to~\cite{MR2209155,MR2373315} for recent results on functionals related to the inequality, that have been obtained by variational methods. These last two papers solve the question, in any dimension, of whether the \emph{first best constant} can be reached. This is equivalent to showing that the difference of the two terms in~\eqref{Ineq:Moser-Trudinger-Onofri} is bounded from below. Earlier results have been obtained by T.~Aubin in~\cite{MR534672}, in the case of the sphere $\S^n$, and P.~Cherrier in~\cite{MR548913} for general $2$-manifolds. The present paper focuses on the value of the \emph{second best constant,} defined as the largest value of $\lambda$ such that~\eqref{Ineq:Moser-Trudinger-Onofri} holds. The value of the \emph{first best constant} is of little concern to us, as it appears as the $\frac14$ coefficient in front of $\nrM{\nabla u}2^2$ and can be factored into $\lambda$. The method used by Z.~Faget in~\cite{MR2209155,MR2373315} relies on a blow-up analysis which is reminiscent for $d=2$ of~\cite{MR1664542}. It generalizes some results contained in~\cite[Theorem 2.50 page 68]{MR681859}. Other references of general interest in the context of the Moser-Trudinger-Onofri inequality are~\cite{MR1143664,MR1230930,MR845999,MR960228,MR3052352}. A review of results related with the Moser-Trudinger-Onofri inequality in the case $\mathcal M=\mathbb S^2$ can be found in~\cite{DEJ}. Let us mention that in~\cite{MR2154301}, A.~Ghigi provides a proof of this inequality based on the Pr\'ekopa-Leindler inequality and that many more details can be found in the book~\cite[Chapters 16-18]{MR3052352} of N.~Ghoussoub and A.~Moradifam. In the context of Einstein-K\"ahler geometry another proof appears in~\cite[Theorem~5.2]{MR2473271} (also see for instance~\cite{Tian:2000fk,phong2008moser} and~\cite{2011arXiv1109.1263B} for recent results on K\"ahlerian manifolds). In~\cite{MR2473271}, Y.A.~Rubinstein gives a proof of the Onofri inequality on $\S^2$ that does not use symmetrization/rearrangement arguments. Also see~\cite{MR2419932} and in particular~\cite[Corollary 10.12]{MR2419932} which contains a reinforced version of the inequality. We shall refer to~\cite{MR2104700} for background material in this direction. The reader interested in understanding how the Moser-Trudinger-Onofri inequality is related to the problem of prescribing the Gaussian curvature on $\S^2$ is invited to refer to~\cite[Section~3]{MR934274} for an introductory survey, and to~\cite{MR925123,MR908146,MR1989228} for more details. There are also many references about the equation~\eqref{Eqn:EL}, which is a particular version of the so-called Liouville equation. The book~\cite{MR2403854} contains a good description of the literature about this equation and many references. More references will be given within the text, whenever needed.

At this point, we should emphasize that in most of the literature the Moser-Trudinger-Onofri inequality in dimension $d=2$ is not written as in~\eqref{Ineq:Moser-Trudinger-Onofri}, but in the form
\[
e^{\,\mu_2\,\nrM{\nabla u}2^2}\ge\,\mathrm C\,\iM{e^u}
\]
for all functions $u\in\mathrm H^1(\mathcal M)$ such that $\iM u=0$, for some constant~$\mathrm C$ which is in general non-explicit. In dimension $d=2$, the optimal constant is $\mu_2=\frac1{4}$. This amounts to write that the functional
\[
u\mapsto\mu_2\,\nrM{\nabla u}2^2+\iM u-\log\(\iM{e^u}\)
\]
is bounded from below by $\log \mathrm C$. The issue of the \emph{first best constant} is to prove that if $\mu_2$ is replaced by a smaller constant, the functional becomes unbounded from below. This is different than proving Inequality~\eqref{Ineq:Moser-Trudinger-Onofri}, except when $\mathrm C=1$ and $\lambda=1$. E.~Onofri proved in~\cite{MR677001} that this is the case, with optimal values for both $\mathrm C$ and $\lambda$, when $\mathcal M=\S^2$, up to a factor $4\,\pi$ that comes from the normalization of $\mathrm{vol}_g(\mathcal M)$. Except for the sphere we are aware of only one occurrence in the literature of the form~\eqref{Ineq:Moser-Trudinger-Onofri} of the inequality, that has been derived by E.~Fontenas in~\cite[Th\'eor\`eme~2]{MR1435336} under more restrictive conditions on $\mathcal M$. This result will be commented in more detail in Remark~\ref{Rem:Fontenas}.

The proof of Theorem~\ref{Thm:Rigidity} is a \emph{rigidity} method inspired by the one of~\cite{MR1134481} for the equation
\[
-\,\Lap u+\lambda\,u=u^{p-1}\,,
\]
which is the Euler-Lagrange equation corresponding to the optimality case in the interpolation inequality 
\be{Ineq:Interp}
\nrM{\nabla v}2^2\ge\frac\lambda{p-2}\,\left[\nrM vp^2-\nrM v2^2\right]\quad\forall\,v\in\H^1(\mathcal M)\,.
\ee
See~\cite{MR1134481,MR1338283,MR1412446,MR2381156,DEKL2014} for further results on this problem and~\cite{MR681859,MR1688256} for general accounts on Sobolev's inequality on Riemannian manifolds. Concerning spectral issues, a standard textbook is~\cite{MR0282313}.

The case of the exponential nonlinearity in~\eqref{Eqn:EL} has been much less considered in the literature, except when $\M$ is the two dimensional sphere $\S^2$. Let us mention the uniqueness result of~\cite{MR1121147} for~\eqref{Eqn:EL} with $\lambda=1$. In~\cite{MR1435336}, and in~\cite{MR1231419} in the case of the ultraspherical operator, the result is achieved by considering the interpolation inequalities~\eqref{Ineq:Interp} and then, as in~\cite{MR1134481} or~\cite{MR1230930} (in the case of the sphere), by taking the limit as \hbox{$p\to\infty$}. Here we consider a direct approach, based on rigidity methods and an associated nonlinear flow. As far as we know, this is an entirely new approach which is interesting in that it provides explicit estimates on the optimal constant in~\eqref{Ineq:Moser-Trudinger-Onofri}.

One may wonder if rigidity results can be achieved for dimensions $d > 2$ with our method. We will give a negative answer in Section~\ref{Sec:Rigidity}. Corollary~\ref{Cor:Rigidity} is established in Section~\ref{Sec:Flow} using a nonlinear flow that has already been considered on the sphere in~\cite{DEJ}. The case $d=1$ is very simple and will be considered for the sake of completeness in Section~\ref{Sec:d=1}. An important application of our method is the case of the Euclidean space with weights, with applications to chemotaxis. Section~\ref{Sec:R2weight} is devoted to this issue with a main result in this direction stated in Theorem~\ref{Thm:Euclidean}, that raises difficult questions of symmetry breaking.

\section{Proof of Theorem~\texorpdfstring{\ref{Thm:Rigidity}}{1}}\label{Sec:Rigidity}

In this section we consider a smooth solution to~\eqref{Eqn:EL} and perform a computation to prove the rigidity result of Theorem~\ref{Thm:Rigidity}. There is no \emph{a priori} reason to assume that $d=2$ and so we shall do the computations for any dimension $d\ge1$, which raises no special additional difficulties. However, due to restrictions that are inherent to the method and will be explicitely exposed, only $d=2$ can be covered. On several occasions, one has to divide by $(d-1)$, so the case $d=1$ has to be excluded and will be handled directly in Section~\ref{Sec:d=1}.

In the case of~\eqref{Ineq:Interp}, it is well known (see~\cite{MR1338283,MR1412446,dolbeault:hal-00784887}) that an interpolation depending on a parameter $\theta\in(0,1)$ between an estimate based on the Ricci curvature and another one based on the first eigenvalue of the Laplace-Beltrami operator can be used to obtain some improvements. Here we apply the same technique and realize in the end that only $\theta=1$ is admissible in dimension $d=2$. However, when $d$ is considered as a real parameter in the range $(1,2)$, it is possible to optimize on $\theta$ when $0\le\theta\le1$. We will comment this and possible improvements at the end of this section.

\subsection*{Preliminaries}
A simple expansion of the square shows that
\[
\|\mathrm H u\|^2=\|\mathrm L u\|^2+\frac1d\,(\Lap u)^2\,.
\]
The Bochner-Lichnerovicz-Weitzenb\"ock formula asserts that 
\[
\frac12\,\Delta\,|\nabla u|^2=\|\mathrm H u\|^2+\nabla (\Lap u)\cdot\nabla u+\Ric(\nabla u,\nabla u)
\]
where $\Ric$ denotes the Ricci tensor and, as a consequence, 
\be{Eqn:BLWIdentity}
\frac12\,\Delta\,|\nabla u|^2=\|\mathrm L u\|^2+\frac1d\,(\Lap u)^2+\nabla (\Lap u)\cdot\nabla u+\Ric(\nabla u,\nabla u)\,.
\ee

\subsection*{An identity based on integrations by parts}
Using integrations by parts, we may notice that
\begin{eqnarray*}
&&-\iM{\Lap u\,|\nabla u|^2\,e^{-u/2}}\\
&&=-\frac12\iM{|\nabla u|^4\,e^{-u/2}}+2\iM{\mathrm H u:\nabla u\otimes\nabla u\,e^{-u/2}}\\
&&=-\frac12\iM{|\nabla u|^4\,e^{-u/2}}+2\iM{\big(\mathrm L u+\frac gd\,\Lap u\big):\nabla u\otimes\nabla u\,e^{-u/2}}\\
&&=-\frac12\iM{|\nabla u|^4\,e^{-u/2}}+2\iM{\mathrm L u:\mathrm M[u]\,e^{-u/2}}\\
&&\hspace*{4cm}+\,\frac2d\iM{\Lap u\,|\nabla u|^2\,e^{-u/2}}\,,
\end{eqnarray*}
with $\mathrm M[u]:=\nabla u\otimes\nabla u-\frac gd\,|\nabla u|^2$, which proves that
\begin{multline*}
\frac{d+2}d\iM{\Lap u\,|\nabla u|^2\,e^{-u/2}}\\
=\frac12\iM{|\nabla u|^4\,e^{-u/2}}-2\iM{\mathrm L u:\mathrm M[u]\,e^{-u/2}}
\end{multline*}
and finally
\begin{multline}\label{Eqn:1}
\iM{\Lap u\,|\nabla u|^2\,e^{-u/2}}\\
=\frac12\,\frac d{d+2}\iM{\!|\nabla u|^4\,e^{-u/2}}-\frac{2\,d}{d+2}\iM{\mathrm L u:\mathrm M[u]\,e^{-u/2}}\,.
\end{multline}

\subsection*{An identity based on the Bochner-Lichnerovicz-Weitzenb\"ock formula}
By expanding $\Lap(e^{-u/2})=(\frac14\,|\nabla u|^2-\frac12\,\Lap u)\,e^{-u/2}$, we have that
\[
\iM{|\nabla u|^2\,\Lap(e^{-u/2})}=\frac14\iM{|\nabla u|^4\,e^{-u/2}}-\frac12\iM{\!\Lap u\,|\nabla u|^2\,e^{-u/2}}
\]
so that, if we multiply~\eqref{Eqn:BLWIdentity} by $e^{-u/2}$ and integrate by parts, then we get
\begin{multline*}
\frac18\iM{|\nabla u|^4\,e^{-u/2}}-\frac14\iM{\Lap u\,|\nabla u|^2\,e^{-u/2}}\\
=\iM{\|\mathrm L u\|^2\,e^{-u/2}}+\frac1d\iM{(\Lap u)^2\,e^{-u/2}}\hspace*{4cm}\\
+\iM{\big(\nabla(\Lap u)\cdot\nabla u\big)\,e^{-u/2}}+\iM{\Ric(\nabla u,\nabla u)\,e^{-u/2}}\\
=\iM{\|\mathrm L u\|^2\,e^{-u/2}}-\frac{d-1}d\iM{(\Lap u)^2\,e^{-u/2}}\\
+\frac12\iM{\Lap u\,|\nabla u|^2\,e^{-u/2}}+\iM{\Ric(\nabla u,\nabla u)\,e^{-u/2}}\,,
\end{multline*}
from which we deduce that
\begin{multline}\label{Eqn:BLW}
\iM{(\Lap u)^2\,e^{-u/2}}\\
=\frac34\,\frac d{d-1}\iM{\Lap u\,|\nabla u|^2\,e^{-u/2}}-\,\frac18\,\frac d{d-1}\iM{|\nabla u|^4\,e^{-u/2}}\\
+\,\frac d{d-1}\iM{\|\mathrm L u\|^2\,e^{-u/2}}+\,\frac d{d-1}\iM{\Ric(\nabla u,\nabla u)\,e^{-u/2}}\,.
\end{multline}

\subsection*{A Poincar\'e inequality}

Since
\[
4\,\Lap(e^{-u/4})=\frac14\,|\nabla u|^2\,e^{-u/4}-\Lap u\,e^{-u/4}\,,
\]
we get that
\begin{multline*}
16\,\iM{|\Lap(e^{-u/4})|^2}=\frac1{16}\iM{|\nabla u|^4\,e^{-u/2}}-\frac12\iM{|\nabla u|^2\,\Lap u\,e^{-u/2}}\\
+\iM{(\Lap u)^2\,e^{-u/2}}\,.
\end{multline*}
On the other hand, a Poincar\'e inequality applied to $\nabla(e^{-u/4})$, as in~\cite[Lemma~7]{dolbeault:hal-00784887}, shows that
\begin{equation}
\label{Ineq:Poincare generic}
\iM{|\Lap(e^{-u/4})|^2}\ge\lambda_1\iM{|\nabla(e^{-u/4})|^2}=\frac{\lambda_1}{16}\iM{|\nabla u|^2\,e^{-u/2}}\,,
\end{equation}
so that
\begin{multline}\label{Ineq:Poincare}
\iM{(\Lap u)^2\,e^{-u/2}}\ge\lambda_1\iM{|\nabla u|^2\,e^{-u/2}}-\frac1{16}\iM{|\nabla u|^4\,e^{-u/2}}\\+\frac12\iM{|\nabla u|^2\,\Lap u\,e^{-u/2}}\,.
\end{multline}

\subsection*{An identity based on the equation}
By expanding $\Lap(e^{-u/2})=(\frac14\,|\nabla u|^2-\frac12\,\Lap u)\,e^{-u/2}$, we have that
\[
\iM{(-\tfrac12\,\Lap u)\,\Lap(e^{-u/2})}=\frac14\iM{(\Lap u)^2\,e^{-u/2}}-\frac18\iM{\Lap u\,|\nabla u|^2\,e^{-u/2}}
\]
so that, if we multiply~\eqref{Eqn:EL} by $\Lap(e^{-u/2})-\frac12\,|\nabla u|^2\,e^{-u/2}$, then we get
\[
\frac14\iM{(\Lap u)^2\,e^{-u/2}}+\frac18\iM{\Lap u\,|\nabla u|^2\,e^{-u/2}}-\frac\lambda2\iM{|\nabla u|^2\,e^{-u/2}}=0\,.
\]
We then split the first term in the equality as the sum of $\frac{1-\theta}4\iM{(\Lap u)^2\,e^{-u/2}}$ and~$\frac\theta4\iM{(\Lap u)^2\,e^{-u/2}}$, with a parameter~$\theta \leq 1$, and we bound the resulting terms using~\eqref{Ineq:Poincare} and~\eqref{Eqn:BLW}, respectively:
\begin{multline*}
\frac{1-\theta}4\left[\lambda_1\iM{|\nabla u|^2\,e^{-u/2}}-\frac1{16}\iM{|\nabla u|^4\,e^{-u/2}}\right.\\
\hspace*{4cm}\left.+\,\frac12\iM{|\nabla u|^2\,\Lap u\,e^{-u/2}}\right]\\
+\frac\theta4\left[\frac34\,\frac d{d-1}\iM{\Lap u\,|\nabla u|^2\,e^{-u/2}}-\,\frac18\,\frac d{d-1}\iM{|\nabla u|^4\,e^{-u/2}}\right.\\
\hspace*{2cm}\left.+\,\frac d{d-1}\iM{\|\mathrm L u\|^2\,e^{-u/2}}+\,\frac d{d-1}\iM{\Ric(\nabla u,\nabla u)\,e^{-u/2}}
\right]\\
+\,\frac18\iM{\Lap u\,|\nabla u|^2\,e^{-u/2}}-\frac\lambda2\iM{|\nabla u|^2\,e^{-u/2}}\le0\,.
\end{multline*}
Collecting terms, we get
\begin{multline*}
\frac\theta4\,\frac d{d-1}\left[\iM{\|\mathrm L u\|^2\,e^{-u/2}}+\iM{\Ric(\nabla u,\nabla u)\,e^{-u/2}}\right]\\
-\,\frac1{64}\(1-\theta+2\,\theta\,\frac d{d-1}\)\iM{|\nabla u|^4\,e^{-u/2}}\hspace*{2cm}\\
+\(\frac{1-\theta}8+\frac{3\,\theta}{16}\,\frac d{d-1}+\frac18\)\iM{\Lap u\,|\nabla u|^2\,e^{-u/2}}\\
+\(\frac{1-\theta}4\,\lambda_1-\frac\lambda2\)\iM{|\nabla u|^2\,e^{-u/2}}\le0
\end{multline*}
and can now use~\eqref{Eqn:1} to obtain
\begin{multline*}
\frac\theta4\,\frac d{d-1}\left[\iM{\|\mathrm L u\|^2\,e^{-u/2}}+\iM{\Ric(\nabla u,\nabla u)\,e^{-u/2}}\right]\\
-\,\frac1{64}\(1-\theta+2\,\theta\,\frac d{d-1}\)\iM{|\nabla u|^4\,e^{-u/2}}\hspace*{2cm}\\
+\(\frac{1-\theta}8+\frac{3\,\theta}{16}\,\frac d{d-1}+\frac18\)\left[\frac12\,\frac d{d+2}\iM{\!|\nabla u|^4\,e^{-u/2}}\right.\\
\left.\hspace*{6cm}-\frac{2\,d}{d+2}\iM{\mathrm L u:\mathrm M[u]\,e^{-u/2}}\right]\\
+\(\frac{1-\theta}4\,\lambda_1-\frac\lambda2\)\iM{|\nabla u|^2\,e^{-u/2}}\le0\,.
\end{multline*}
Recall that $\mathrm M[u]$ denotes the trace free tensor
\[
\mathrm M[u]:=\nabla u\otimes\nabla u-\frac gd\,|\nabla u|^2\,.
\]
We observe that
\[\label{Eqn:nabla4}
\|\mathrm M[u]\|^2=\left\|\,\nabla u\otimes\nabla u-\frac gd\,|\nabla u|^2\,\right\|^2=\(1-\frac1d\)|\nabla u|^4\,. 
\]
Altogether we end up with
\begin{multline*}
\iM{\Big(\mathsf a\,\|\mathrm L u\|^2+\mathsf b\,(\mathrm L u:\mathrm M[u])+\mathsf c\,\|\mathrm M[u]\|^2\Big)\,e^{-u/2}}\\
+\frac\theta4\,\frac d{d-1}\iM{\Ric(\nabla u,\nabla u)\,e^{-u/2}}\\
+\(\frac{1-\theta}4\,\lambda_1-\frac\lambda2\)\iM{|\nabla u|^2\,e^{-u/2}}\le0
\end{multline*}
with
\begin{eqnarray*}
&&\mathsf a=\frac\theta4\,\frac d{d-1}\,,\\
&&\mathsf b=-\(\frac{1-\theta}8+\frac{3\,\theta}{16}\,\frac d{d-1}+\frac18\)\frac{2\,d}{d+2}\,,\\
&&\mathsf c=\left[\(\frac{1-\theta}8+\frac{3\,\theta}{16}\,\frac d{d-1}+\frac18\)\frac12\,\frac d{d+2}-\,\frac1{64}\(1-\theta+2\,\theta\,\frac d{d-1}\)\right]\,\frac d{d-1}\,.
\end{eqnarray*}

\begin{Rem}\label{Rem:Lichnerowicz}
By the Lichnerowicz' theorem (see~\cite{MR0124009} or~\cite[Section~2]{dolbeault:hal-00784887})
\[
\frac d{d-1}\iM{\Ric(\nabla u,\nabla u)\,e^{-u/2}}\le\iM{|\nabla u|^2\,e^{-u/2}}
\]
so that the largest possible value of $\lambda$ for which we \emph{a priori} know that
\[
\frac\theta4\,\frac d{d-1}\iM{\Ric(\nabla u,\nabla u)\,e^{-u/2}}+\(\tfrac 14\,\lambda_1\,(1-\theta)-\lambda\)\iM{|\nabla u|^2\,e^{-u/2}}
\]
is nonnegative corresponds to the smallest possible value of $\theta$, \emph{i.e.}~$\theta=\theta_0(d)$.
\end{Rem}

\subsection*{Discussion}
With a simple but tedious computation, one can show that the discriminant $\delta:=\mathsf b^2-\,4\,\mathsf a\,\mathsf c$ has the sign of
\[
16\,(d-1)^2-(6-d)\,(d+2)\,\theta\,.
\]
If we denote by $\theta_0=\theta_0(d)$ the value of $\theta$ for which $\delta=0$, then we have 
\[
\theta_0=\frac{16\,(d-1)^2}{(6-d)\,(d+2)}\,.
\]
Altogether we can rewrite our estimate as
\begin{multline*}
\mathsf a\iM{\left\|\,\mathrm L u + \tfrac{\mathsf b}{2\,\mathsf a}\,\mathrm M[u]\,\right\|^2\,e^{-u/2}}+\(\mathsf c-\tfrac{\mathsf b^2}{4\,\mathsf a}\)\iM{\|\mathrm M[u]\|^2\,e^{-u/2}}\\
+\,4\,(1-\theta)\iM{|\Lap(e^{-u/4})|^2}+\frac\theta4 \,\frac d{d-1}\iM{\Ric(\nabla u,\nabla u)\,e^{-u/2}}\\
-\frac{\lambda}2 \iM{|\nabla u|^2\,e^{-u/2}}=0
\end{multline*}
and use the Poincar\'e inequality~\eqref{Ineq:Poincare generic} to establish that
\begin{multline}\label{estttt}
0\ge\mathsf a\iM{\left\|\,\mathrm L u+\tfrac{\mathsf b}{2\,\mathsf a}\,\mathrm M[u]\,\right\|^2\,e^{-u/2}}-\tfrac\delta{4\,\mathsf a}\iM{\|\mathrm M[u]\|^2\,e^{-u/2}}\\
+\(\tfrac14\,\lambda_1\,(1-\theta)-\tfrac\lambda2\)\iM{|\nabla u|^2\,e^{-u/2}}+\tfrac\theta4\,\tfrac d{d-1}\iM{\Ric(\nabla u,\nabla u)\,e^{-u/2}}\,.
\end{multline}
Our goal is to show that $u$ has to be a constant, that is, $\iM{|\nabla u|^2\,e^{-u/2}}=0$. We assume that $d\ge2$ is an integer. The discriminant $\delta$ is nonpositive if and only if $d<6$ and $\theta\ge\theta_0(d)$. This is compatible with the condition $\theta\le1$ only if $d=2$; in that case, $\theta=\theta_0(2)=1$ and as a consequence, we can rewrite~\eqref{estttt} as
\begin{multline*}
0\ge \iM{\|\mathrm L u-\tfrac12\,\mathrm M[u]\|^2\,e^{-u/2}}+\iM{\Ric(\nabla u,\nabla u)\,e^{-u/2}}\\
-\lambda\iM{|\nabla u|^2\,e^{-u/2}}\\
\ge(\lambda_\star-\lambda)\iM{|\nabla u|^2\,e^{-u/2}}\,.
\end{multline*}
Hence we have shown that $\nabla u\equiv0$ for any $\lambda<\lambda_\star$, which concludes the proof of Theorem~\ref{Thm:Rigidity}.

\begin{Rem}\label{Rem:Fontenas} In order to compare our results with the results deduced from the \emph{curvature-dimension} method (see for instance~\cite{MR3155209}), we can consider the case where $\Ric$ is uniformly bounded from below by some positive constant~$\rho$ and formally assume that $d\in (1,2)$ takes real values. This can be made precise for instance in the setting of the ultra-spherical operator (see for instance~\cite{MR1231419}), with exactly the same conditions as above. See~\cite{MR1435336} and~\cite[Section~7.1]{DEJ} for more details. If $1<d\le2$, we find that rigidity holds if 
\[
\lambda\le\max_{\theta\in[\theta_0(d),1]}\(\frac12\,\lambda_1\,(1-\theta)+\frac{\theta}2\,\frac d{d-1}\,\rho\)=\frac12\,\lambda_1\,(1-\theta_0(d))+\frac{\theta_0(d)}2\,\frac d{d-1}\,\rho
\]
according to Remark~\ref{Rem:Lichnerowicz}. Let $x=\frac d{d-1}\,\frac\rho{\lambda_1}\le1$. We have found that that rigidity holds~if
\[
2\,\frac\lambda{\lambda_1}\le1-\theta_0(d)+\,\theta_0(d)\,x=:f_1(x)
\]
Quite surprisingly, a better condition has been obtained in~\cite[Th\'eor\`eme~2]{MR1435336}, when $1<d<2$, which amounts to
\[
2\,\frac\lambda{\lambda_1}\le d\,(2-d)+(d-1)^2\,x=:f_2(x)\,,
\]
by taking the limit as $p\to\infty$ in~\eqref{Ineq:Interp}. We may indeed check that $f_2(x)-f_1(x)=\frac{(d-1)^2\,(d-2)^2}{(6-d)\,(d+2)}\,(1-x)\ge0$

Without assuming the positivity of $\rho$, one gets a similar result with our approach. In the range $d\in(1,2)$, our computations show that rigidity holds for any~$\lambda$ at most equal to the infimum on $u\in\H^2(\M)\setminus\{0\}$ of 
\[
2 \iM{\Big[\,\mathsf a\,\left\|\,\mathrm L u+\tfrac{\mathsf b}{2\,\mathsf a}\,\mathrm M[u]\,\right\|^2+\,{\scriptstyle 4\,(1-\theta)}\,|\Lap(e^{-u/4})|^2+\tfrac\theta4\,\tfrac d{d-1}\,\Ric(\nabla u,\nabla u)\Big]\,e^{-u/2}}
\]
under the condition that $\iM{|\nabla u|^2\,e^{-u/2}}=1$, $\mathsf a=\frac{4\,d\,(d-1)}{(6-d)\,(d+2)}$, $\mathsf b= - \frac{d\,(3\,d+2)}{2\,(6-d)\,(d+2)}$ and $\theta=\theta_0(d)$. However, in the same spirit as above, a passage to the limit as \hbox{$p\to\infty$} in the inequality obtained in~\cite[Theorem~4]{dolbeault:hal-00784887} gives a better result. 

Let us emphasize that these considerations are essentially formal because $d$ is restricted to the interval $(1,2)$ but can be entirely justified, as it is currently done in the \emph{curvature-dimension} approach. See for instance~\cite[Th\'eor\`eme~2]{MR1435336}, and related references.\end{Rem}

\section{Proof of Corollary~\texorpdfstring{\ref{Cor:Rigidity}}{2}}\label{Sec:Flow}

Let us minimize the functional
\[
\mathcal F_\lambda[u]:=\frac14\,\nrM{\nabla u}2^2+\lambda\,\iM u-\lambda\,\log\(\iM{e^u}\)
\]
on $\H^1(\M)$. According to~\cite{MR2209155,MR2373315}, $\mathcal F_\lambda[u]$ has no minimizer if $\lambda>1$. Let us assume that
\[
\lambda<1
\]
(we shall take care of the equality case later). It is then standard that there is a well-defined minimizer $u$. Note that since $\mathcal F_\lambda[u] $ does not change when adding some constant to $u$, we can choose the value of $\iM{e^u} = \lambda$ and then $u$ satisfies~\eqref{Eqn:EL}. If $u$ is smooth and $\lambda<\lambda_\star$, we can apply the result of Theorem~\ref{Thm:Rigidity}. Then the minimizer $u$ has to be a constant, for instance $u\equiv1$, so that 
\[
\mathcal F_\lambda[u]\ge\mathcal F_\lambda[1]=0\quad\forall\,u\in\H^1(\M)\,.
\]
 Notice that we can get rid of any smoothness requirement by considering the flow below. By passing to the limit as $\lambda\nearrow\lambda_\star$, we get that the inequality also holds true if $\lambda=\lambda_\star$.

Using the following Taylor expansion of~$\mathcal F_\lambda$ as $\varepsilon\to0$,
\begin{align*}
\mathcal F_\lambda[1+\varepsilon\,\varphi]&=\varepsilon^2\left[\frac14\,\nrM{\nabla\varphi}2^2+\frac\lambda2\,\iM{\varphi^2}\right]+o(\varepsilon^2)
\\
&= \frac{\varepsilon^2}{2} \left(\frac{\lambda_1}{2} - \lambda\right) +o(\varepsilon^2)\,,
\end{align*}
where $\varphi$ is an eigenfunction associated with the first positive eigenvalue~$\lambda_1$ of $-\,\Lap$, it is straightforward to see that the best constant $\lambda$ in~\eqref{Ineq:Moser-Trudinger-Onofri} is such that
\[
\lambda\le\frac{\lambda_1}2\,.
\]
To complete the proof of Corollary~\ref{Cor:Rigidity}, it remains to consider the case $\lambda_\star<\lambda_1/2$ and show that the optimal constant cannot be equal to $\lambda_\star$. This can be done in the same spirit as in~\cite[Corollary~2]{dolbeault:hal-00784887}. Let us consider the evolution equation defined by

\be{Eqn:Evol}
\frac{\partial f}{\partial t}=\Lap(e^{-f/2})-\tfrac12\,|\nabla f|^2\,e^{-f/2}\;,
\ee
with initial datum $u\in\H^1(\mathcal M)$. Let us define
\begin{multline*}
\mathcal G_\lambda[f]:=\iM{\|\,\mathrm L f-\tfrac12\,\mathrm M[f]\,\|^2\,e^{-f/2}}+\iM{\Ric(\nabla f,\nabla f)\,e^{-f/2}}\\
-\lambda\iM{|\nabla f|^2\,e^{-f/2}}\,.
\end{multline*}
Then for any $\lambda\le\lambda_\star$ we have
\[
\frac d{dt}\mathcal F_\lambda[f(t,\cdot)]=\iM{\(-\tfrac12\,\Lap f+\lambda\)\(\Lap(e^{-f/2})-\tfrac12\,|\nabla f|^2\,e^{-f/2}\)}=-\,\mathcal G_\lambda[f(t,\cdot)]
\]
Since $\mathcal F_\lambda$ is nonnegative and $\lim_{t\to\infty}\mathcal F_\lambda[f(t,\cdot)]=0$, we obtain that
\[
\mathcal F_\lambda[u]\ge\int_0^\infty\mathcal G_\lambda[f(t,\cdot)]\,dt
\]
for any solution $f$ to~\eqref{Eqn:Evol} with initial datum $u\in\L^1(\mathcal M)$ is such that $\nabla u\in\L^2(\mathcal M)$. We have an equality if the solution is smooth for any $t\ge0$. Otherwise we have to regularize and then pass to the limit so that, with full generality, we can only expect for an inequality.

\begin{Rem}\label{rem:positivity lambda}One has to mention that the sphere $\mathcal M=\S^2$ is an important case of application of our method, for which other types of remainder terms can be produced. See~\cite{DEJ} for more details. It has to be noted that on $\S^2$ we have \hbox{$\lambda_\star=\rho=\lambda_1/2=1$}. As noted by many authors, \emph{e.g.}, in~\cite{MR1230930,MR1231419,MR1435336,DEJ} (also see references in~\cite{DEJ}), the Onofri inequality is a limit case of various Sobolev type inequalities, for which similar methods have been developed: see~\cite{DEKL2014} for a review and some recent results. 

Another interesting case for our method is the flat square torus, defined as the square $(0,1)\times (0,1)$ with double periodicity, or simply $\mathbb T^2=\R^2/\mathbb Z^2$. For general properties of the solutions of~\eqref{Eqn:EL} on $\mathbb T^2$, see~\cite[Section~2.5]{MR2403854}. In this case, and also for any compact manifold $\M$ with nonnegative Ricci curvature, by defining $v=e^{-u/2}$ one can see that
\[\label{LambdaStartorus}
\lambda_\star=\inf_{v\in\H^2(\M)\setminus\{0\}}\frac{\iM{\big(\frac{\,\|\,\mathrm L v\,\|^2}v+\frac1v\,\Ric(\nabla v,\nabla v)\big)}}{\iM{\frac{|\nabla v|^2}{v}}}\,.
\]
According to~\cite[Theorem~0.1]{MR1673972}, there exists~$C>0$ such that for any $0<\lambda<C$, all solutions $u$ of~\eqref{Eqn:EL} are uniformly bounded in $L^\infty(\M)$. In the proofs of Theorem~\ref{Thm:Rigidity} and Corollary~\ref{Cor:Rigidity}, it is henceforth possible to replace $\lambda_\star$ by the infimum taken over the space of all functions $v\in\H^2(\M)\setminus\{0\}$ such that $-\,2\,\mathrm{var}(\log v)=\mathrm{var}(u)=\mathrm{supess}(u)-\mathrm{infess}(u)\le\mathcal K$ for some positive constant $\mathcal K$, which is independent of $\lambda>0$. Then, with this new definition of $\lambda_\star$, we obtain the estimate
\[
\lambda_\star\ge e^{-\mathcal K/2}\,\inf\frac{\iM{\big(\|\,\mathrm L v\,\|^2+\Ric(\nabla v,\nabla v)\big)}}{\iM{|\nabla v|^2}}=e^{-\mathcal K/2}\,\inf\frac{\iM{\big(\Delta v)^2}}{\iM{|\nabla v|^2}}\,,
\]
where the last equality is a straightforward consequence of the Bochner-Lichne\-rovicz-Weitzenb\"ock formula. As a consequence of the Poincar\'e inequality (see for instance~\cite[Lemma~5]{dolbeault:hal-00784887}), we obtain the estimate $\lambda_\star\ge e^{-\mathcal K/2}\,\lambda_1$. It is remarkable that our method applies when the lowest eigenvalue of $\Ric(\nabla v,\nabla v)$ takes value $0$ in $\M$. In the case of the flat torus, we even have that $\Ric(\nabla v,\nabla v)\equiv0$ on~$\M$.\end{Rem}

\section{The case \texorpdfstring{$d=1$}{d=1}}\label{Sec:d=1}

For simplicity, we can consider the case of the circle. Hence we identify $\mathcal M$ with the $1$-periodic interval \hbox{$[0,1)\approx\R/\mathbb Z\approx\S^1$}. Consider a
solution of the ordinary differential equation
\be{Eqn:1d}
-\,\frac12\,u''+\lambda-e^u=0
\ee
with periodic boundary conditions. If we multiply the equation by $(e^{-u/2})''-\frac12\,|u'|^2\,e^{-u/2}$, then we get that
\[
\int_0^1\(\tfrac14\,|u''|^2+\tfrac18\,|u'|^2\,u''-\,\tfrac\lambda2\,|u'|^2\)e^{-u/2}\,dx=0\,.
\]
The middle term is easy to handle using one integration by parts:
\be{intparts}
\int_0^1|u'|^2\,u''\,e^{-u/2}\,dx=\tfrac16\int_0^1|u'|^4\,e^{-u/2}\,dx\,.
\ee
Hence we have
\be{ineq48}
\int_0^1\(\tfrac14\,|u''|^2+\tfrac1{48}\,|u'|^4-\,\tfrac\lambda2\,|u'|^2\)e^{-u/2}\,dx=0\,.
\ee
On the other hand, by the Poincar\'e inequality, 
\[
\int_0^1\left|\(e^{-u/4}\)''\right|^2\,dx \ge 4\,\pi^2\,\int_0^1\left|\(e^{-u/4}\)'\right|^2\,dx\,,
\]
where $4\,\pi^2$ is the first positive eigenvalue of $-\frac{d^2}{dx^2}$ on the periodic interval of length~$1$.
From~\eqref{intparts} we derive
\be{spectralineq}
\int_0^1\(|u''|^2-\tfrac1{48}\,|u'|^4\)\,e^{-u/2} \,dx - 4\,\pi^2\int_0^1|u'|^2\,e^{-u/2} \,dx \ge 0\,.
\ee
Combining~\eqref{ineq48} and~\eqref{spectralineq}, we get 
\[
\frac{5}{96}\,\int_0^1|u'|^4\,e^{-u/2}\,dx+(2\pi^2-\lambda)\,\int_0^1\,|u'|^2\,e^{-u/2}\,dx \le 0\,.
\]
Hence we have proven the following result.
\begin{Prop}\label{Prop:Rigidity1d} Assume that $d=1$. With the above notations, if $u$ is a smooth solution to~\eqref{Eqn:1d} on $\S^1\approx[0,1)$, then $u$ is a constant function for any $\lambda\in(0,2\,\pi^2)$. \end{Prop}
Exactly as in the case of a manifold of dimension two, a variational approach allows to deduce a Moser-Trudinger-Onofri inequality.
\begin{Cor}\label{Cor:Rigidity1d} If $d=1$, then the following inequality holds on $\S^1\approx[0,1)$:
\[
\frac1{8\,\pi^2}\,\int_0^1|u'|^2\,dx+\int_0^1u\,dx\ge\,\log\(\int_0^1e^u\,dx\)\quad\forall\,u\in\H^1(\S^1)\,.
\]
Moreover $8\,\pi^2$ is the optimal constant.\end{Cor}
The only difference with Corollary~\ref{Cor:Rigidity} is that we can identify the optimal constant in the inequality by considering $u=1+\varepsilon\,\varphi$ and by taking the limit as $\varepsilon\to0$, with $\varphi(x)=\cos(\,2\pi\,x)$.

\section{Weighted Moser-Trudinger-Onofri inequalities on the two-dimensional Euclidean space}\label{Sec:R2weight}

The Euclidean Onofri inequality (see~\cite{MR1143664,DEJ}) can be deduced from~\eqref{Ineq:Moser-Trudinger-Onofri} when $\mathcal M=\S^2$ using the stereographic projection and reads
\be{Onofri:Euclidean}
\frac1{16\,\pi}\irdeux{|\nabla u|^2}\ge\log\(\irdmu{e^u}\)-\irdmu u\,.
\ee
Here $d\mu(x)\,dx$ denotes the probability measure defined by $\mu(x)=\frac1\pi\,(1+|x|^2)^{-2}$, $x\in\R^2$, and the inequality holds for any function $u\in\L^1(\R^2,d\mu)$ such that $\nabla u\in\L^2(\R^2)$. The constant $16\,\pi$ is optimal as can be shown by considering the inequality on $\S^2$ and comparing with the value given when expanding around a constant, as was done in Section~\ref{Sec:Flow}.

\medskip In this section, our goal is to give sufficient conditions on a general probability measure $\mu$ so that the inequality
\be{eq:weighted_onofri}
\frac1{16\,\pi}\irdeux{|\nabla u|^2}\ge\lambda\left[\log\(\irdmu{e^u}\)-\irdmu u\right]
\ee
holds for some $\lambda>0$ and get an estimate of the optimal value of $\lambda$. Here $d\mu=\mu\,dx$ is a probability measure with density $\mu$ with respect to the Lebesgue measure. All our computations are done without symmetry assumption, and our final estimate is~\eqref{Eqno:EuclideanCondition}. In practical applications (see Examples~\ref{Ex1}--\ref{Ex4}) the function $\mu$ is radially symmetric \emph{and} one has to assume that $\lambda$ is in a range for which the solution to~\eqref{eq:weighted_onofri}, or at least the optimal function for~\eqref{Onofri:Euclidean}, is radially symmetric. This delicate issue of \emph{symmetry breaking} will be illustrated in Example~\ref{Ex3}.

\medskip Since~\eqref{eq:weighted_onofri} does not change when adding some constant to $u$, we can look for minimizers satisfying the constraint $\irdmu{e^u} = 1$.
These solve the Euler-Lagrange equation
\begin{equation}
-\frac1{8\,\pi}\,\Delta u + \lambda\,\mu - \lambda\,e^u\,\mu = 0\,.
\label{eq:EL_weigthed}
\end{equation}

We can multiply each term of~\eqref{eq:EL_weigthed} by $ \frac 1\mu\,\Delta ( e^{-\frac u2} ) - \frac 1{2\,\mu} \left|\nabla u\right|^2 e^{-\frac u2}$ and integrate,
which gives the following identities
\[
-\irdeux{ \Delta u\, \Delta (e^{-u/2})\,\frac 1\mu} = \irdeux{ \Delta u \left( \Delta u - \frac 12 \left|\nabla u\right|^2 \right) e^{-u/2}\,\frac 1{2\,\mu}}\,,
\]
\begin{gather*}
\irdeux{ \mu\,\frac 1\mu\,\Delta(e^{-u/2})} = 0\,,
\\
\irdeux{ e^u \mu \(\frac 1\mu\,\Delta ( e^{-\frac u2} ) - \frac 1{2\,\mu} \left|\nabla u\right|^2 e^{-\frac u2}\)} = 0\,.
\end{gather*}
Defining $\nu := {e^{-u/2}}/{\mu} = e^{-u/2-g}$ with $g:=\log\mu$ and $d\nu := \nu\,dx$ we have
\[
\mathcal I[u] = 2\irdnu{ (\Delta u)^2} +\irdnu{ \Delta u \left| \nabla u \right|^2}-16\,\pi\,\lambda\irdeux{ \left|\nabla u \right|^2 e^{-u/2}} = 0\,.
\]
Let us introduce some notations, which are consistent with the ones on manifolds. Let us denote by $\mathrm Hu=\big(\frac{\partial^2u}{\partial x_i\partial x_j}\big)_{i,j=1,2}$ the Hessian of $u$, $\mathrm Lu=\mathrm Hu-\frac12\,\Delta u\,\mathrm I_2$ is the trace free Hessian and $\mathrm M[u]:=\nabla u \otimes \nabla u - \frac 12\, \left|\nabla u \right|^2\,\mathrm I_2$, where $\nabla u \otimes \nabla u=\big(\frac{\partial u}{\partial x_i}\,\frac{\partial u}{\partial x_j}\big)_{i,j=1,2}$. For the convenience of the reader, we split the computations in four steps.

\medskip\noindent 1) Let us start with some preliminary computations. An integration by parts shows that
\begin{multline}\label{Eqn:Elim1}
2\irdnu{\Delta u\,\nabla u\cdot\nabla g}-\irdnu{\left|\nabla u\right|^2\,\nabla u\cdot\nabla g}\\
\hspace*{-1.5cm}=-\,2\irdnu{\mathrm Hu:(\nabla u\otimes\nabla g)}-\,2\irdnu{\(\mathrm Hg-\nabla g\otimes\nabla g\):(\nabla u\otimes\nabla u)}\\
=-\,2\irdnu{\mathrm Hu:(\nabla u\otimes\nabla g)}-\,2\irdnu{\mathrm Hg:(\nabla u\otimes\nabla u)}+2\irdnu{(\nabla u\cdot\nabla g)^2}\,.
\end{multline}
By expanding $\mathrm Lu - \frac 12\,\mathrm M[u]$, we also get that
\begin{multline}
\label{Eqn:Elim2}
\irdnu{\(\mathrm Lu - \frac 12\,\mathrm M[u]\):(\nabla u\otimes\nabla g)}\\
=\irdnu{\mathrm Hu:(\nabla u\otimes\nabla g)}-\frac12\irdnu{\Delta u\,\nabla u\cdot\nabla g} - \frac 14\,\irdnu{\left|\nabla u\right|^2\,\nabla u\cdot\nabla g}\,.
\end{multline}
Recalling the definition of $\nu=e^{-\frac12u-g}$, we also find that
\begin{multline}\label{Eqn:Elim3}
-\frac12\irdnu{\left|\nabla u\right|^2\,\nabla u\cdot\nabla g}=\irdeux{\left|\nabla u\right|^2\,\nabla g\,e^{-g}\cdot\nabla(e^{-\frac12u})}\\
=-\,2\irdnu{\mathrm Hu:(\nabla u\otimes\nabla g)}+\irdnu{\left|\nabla u\right|^2\,(|\nabla g|^2-\Delta g)}\,.
\end{multline}
Equations~\eqref{Eqn:Elim1},~\eqref{Eqn:Elim2} and~\eqref{Eqn:Elim3} allow us to eliminate
\[
\irdnu{\Delta u\,\nabla u\cdot\nabla g}\,,\quad\irdnu{\left|\nabla u\right|^2\,\nabla u\cdot\nabla g}\quad\mbox{and}\quad\irdnu{\mathrm Hu:(\nabla u\otimes\nabla g)}
\]
in terms of the other quantities:
\begin{eqnarray}
&&\irdnu{\Delta u\,\nabla u\cdot\nabla g}=\irdnu{\left|\nabla u\right|^2\,(|\nabla g|^2-\Delta g)}\nonumber\\
&&\hspace*{4cm}-\,2\irdnu{\(\mathrm Lu - \frac 12\,\mathrm M[u]\):(\nabla u\otimes\nabla g)}\,,\label{a}\\
&&\irdnu{\left|\nabla u\right|^2\,\nabla u\cdot\nabla g}=4\irdnu{\mathrm Hg:(\nabla u\otimes\nabla u)}-\,4\irdnu{(\nabla u\cdot\nabla g)^2}\nonumber\\
&&\hspace*{4cm}-\,8\irdnu{\(\mathrm Lu - \frac 12\,\mathrm M[u]\):(\nabla u\otimes\nabla g)}\label{b}\\
&&\hspace*{4cm}+\,6\irdnu{\left|\nabla u\right|^2\,(|\nabla g|^2-\Delta g)}\,,\nonumber\\
&&\irdnu{\mathrm Hu:(\nabla u\otimes\nabla g)}=\irdnu{\mathrm Hg:(\nabla u\otimes\nabla u)}-\irdnu{(\nabla u\cdot\nabla g)^2}\nonumber\\
&&\hspace*{4cm}-\,2\irdnu{\(\mathrm Lu - \frac 12\,\mathrm M[u]\):(\nabla u\otimes\nabla g)}\label{c}\\
&&\hspace*{4cm}+\,2\irdnu{\left|\nabla u\right|^2\,(|\nabla g|^2-\Delta g)}\,.\nonumber
\end{eqnarray}
Moreover, the reader is invited to check that
\be{d}
\|\mathrm Lu\|^2=\|\mathrm Hu\|^2-\frac 12\,(\Delta u)^2
\ee
and
\be{e}
\|\mathrm M[u]\|^2=\frac 12\,|\nabla u|^4\,.
\ee

\medskip\noindent 2) On the one hand, integrating the second term in the expression of $\mathcal I[u]$ by parts gives
\begin{align*}
&\irdnu{ \Delta u \left| \nabla u \right|^2}\\
&= - \irdeux{ \nabla u \cdot \nabla \left( \left| \nabla u \right|^2 \nu \right)}
\\
&= - \irdnu{ \nabla u\cdot\left(
2\,\mathrm Hu\,\nabla u
- \frac12\,\nabla u \left|\nabla u\right|^2
- \left| \nabla u \right|^2 \nabla g
\right)}
\\
&= - \irdnu{ \left(
2\,\mathrm Hu : \nabla u \otimes \nabla u
- \frac12 \left|\nabla u\right|^4
- \left| \nabla u \right|^2\,(\nabla u \cdot\nabla g)
\right)}
\\
&= - \irdnu{ \left(
2\,\mathrm Lu : \nabla u \otimes \nabla u
+ \Delta u \left| \nabla u\right|^2
- \frac12 \left|\nabla u\right|^4
- \left| \nabla u \right|^2\,(\nabla u \cdot\nabla g)
\right)}
\\
&= - \irdnu{ \left(
2\,\mathrm Lu : \mathrm M[u]
+ \Delta u \left| \nabla u\right|^2
- \frac12 \left|\nabla u\right|^4 
- \left| \nabla u \right|^2\,(\nabla u \cdot\nabla g)\)}\,,
\end{align*}
that is
\[
\irdnu{ \Delta u \left| \nabla u \right|^2}
= - \irdnu{ \left(
\mathrm Lu : \mathrm M[u]
- \frac14 \left|\nabla u\right|^4
- \frac12 \left| \nabla u \right|^2\,(\nabla u \cdot\nabla g)
\)}\,.
\]
According to~\eqref{b} and~\eqref{e}, we obtain
\begin{multline}
\label{eq:ipp_1}
\irdnu{\Delta u \left| \nabla u \right|^2}=-\irdnu{\mathrm Lu : \mathrm M[u]}\\
\hspace*{1cm}+\frac12\irdnu{\|\mathrm M[u]\|^2}+\,2\irdnu{\mathrm Hg:(\nabla u\otimes\nabla u)}\\
\hspace*{3.5cm}-\,2\irdnu{(\nabla u\cdot\nabla g)^2}-\,4\irdnu{\(\mathrm Lu - \frac 12\,\mathrm M[u]\):(\nabla u\otimes\nabla g)}\\+\,3\irdnu{\left|\nabla u\right|^2\,(|\nabla g|^2-\Delta g)}\,.
\end{multline}

\medskip\noindent 3) On the other hand, integrating by parts twice yields
\begin{align}
\label{eq:ipp_2}
\irdnu{ \Delta \left| \nabla u \right|^2}
&= \int_{\R^2}\left|\nabla u\right|^2 \Delta (e^{-u/2})\,\frac {dx}{\mu}
+ \irdnu{ \left|\nabla u\right|^2 \,\(|\nabla g|^2-\Delta g\)}
\nonumber
\\
&\hspace*{1cm}+ \irdnu{ \left|\nabla u\right|^2\,(\nabla u \cdot\nabla g)}\,,
\nonumber
\\
&= -\frac 12\irdnu{ \Delta u \left|\nabla u\right|^2}+ \frac 14 \irdnu{ \left|\nabla u\right|^4}
\nonumber
\\
&\hspace*{1cm}
+ \irdnu{ \left|\nabla u\right|^2 \,\(|\nabla g|^2-\Delta g\)}
+ \irdnu{ \left|\nabla u\right|^2\,(\nabla u \cdot\nabla g)}\,.
\end{align}
Integrating by parts again we have that
\[
\irdnu{ \nabla \Delta u \cdot \nabla u} =
- \irdnu{ \Delta u \( \Delta u - \frac 12 \left|\nabla u\right|^2 - (\nabla u \cdot\nabla g) \)}\,,
\]
which we can use along with the Bochner-Lichnerovicz-Weitzenb\"ock formula on~$\mathbb \R^2$ (with Ricci tensor identically equal to $0$),
\[
\Delta \left| \nabla u \right|^2 = 2\,\| \mathrm Lu \|^2 + (\Delta u)^2 + 2\, \nabla \Delta u\cdot\nabla u \,,
\]
to get
\begin{align*}
\irdnu{ \Delta \left| \nabla u \right|^2}
&= \irdnu{\( 2\,\| \mathrm Lu \|^2 - (\Delta u)^2 + \Delta u \left|\nabla u\right|^2 + 2\,\Delta u\,(\nabla u \cdot\nabla g)\)}\,.
\end{align*}
Combined with~\eqref{eq:ipp_2} this proves that
\begin{multline*}
\irdnu{ (\Delta u)^2 } = 2 \irdnu{ \| \mathrm Lu \|^2 } + \frac 32 \irdnu{ \Delta u \left|\nabla u\right|^2 }
- \frac 14 \irdnu{ \left| \nabla u \right|^4 }
\\
- \irdnu{ \left| \nabla u\right|^2 \,\(|\nabla g|^2-\Delta g\)}
+\, 2 \irdnu{ \Delta u\,(\nabla u \cdot\nabla g) }
- \irdnu{ \left| \nabla u\right|^2\,(\nabla u \cdot\nabla g) }\,.
\end{multline*}
Using~\eqref{a},~\eqref{b},~\eqref{e} and~\eqref{eq:ipp_1}, we obtain
\begin{eqnarray*}
&&\hspace*{-1cm}\irdnu{ (\Delta u)^2 }\\ &=& 2 \irdnu{ \| \mathrm Lu \|^2 } 
\\
&&-\,\frac32\irdnu{\mathrm Lu : \mathrm M[u]}+\frac34\irdnu{\|\mathrm M[u]\|^2}+\,3\irdnu{\mathrm Hg:(\nabla u\otimes\nabla u)} \\
&&\hspace*{2cm}-\,3\irdnu{(\nabla u\cdot\nabla g)^2}-\,6\irdnu{\(\mathrm Lu - \frac 12\,\mathrm M[u]\):(\nabla u\otimes\nabla g)}
\\
&&\hspace*{2cm}+\,\frac92\irdnu{\left|\nabla u\right|^2\,(|\nabla g|^2-\Delta g)}\\
&&-\,\frac 12 \irdnu{\|\mathrm M[u]\|^2}
- \irdnu{ \left| \nabla u\right|^2 \,\(|\nabla g|^2-\Delta g\)}
\\
&&
+\,2\irdnu{\left|\nabla u\right|^2\,(|\nabla g|^2-\Delta g)}-\,4\irdnu{\(\mathrm Lu - \frac 12\,\mathrm M[u]\):(\nabla u\otimes\nabla g)}
\\
&&
-\,4\irdnu{\mathrm Hg:(\nabla u\otimes\nabla u)}+\,4\irdnu{(\nabla u\cdot\nabla g)^2} \\
&&\hspace*{2cm}+\,8\irdnu{\(\mathrm Lu - \frac 12\,\mathrm M[u]\):(\nabla u\otimes\nabla g)} \\
&&\hspace*{2cm}-\,6\irdnu{ \left| \nabla u\right|^2 \,\(|\nabla g|^2-\Delta g\)}\,.
\end{eqnarray*}
Collecting terms, we arrive at
\begin{multline}\label{Eqn:Deltau2}
\irdnu{(\Delta u)^2}=2\irdnu{\|\mathrm Lu\|^2}-\,\frac32\irdnu{\mathrm Lu:\mathrm M[u]}+\,\frac14\irdnu{\|\mathrm M[u]\|^2}\\-\,2\irdnu{\(\mathrm Lu-\frac12\,\mathrm M[u]\):(\nabla u\otimes\nabla g)}\\
\hspace*{2cm}-\irdnu{\mathrm Hg:(\nabla u\otimes\nabla u)}+\irdnu{(\nabla u\cdot\nabla g)^2}\\
-\,\frac12\irdnu{\left|\nabla u\right|^2\,(|\nabla g|^2-\Delta g)}\,.
\end{multline}

\medskip\noindent 4) By reinjecting~\eqref{eq:ipp_1} and~\eqref{Eqn:Deltau2} in the expression of $\mathcal I$, we get that
\begin{multline*}
0=\mathcal I[u]=4\irdnu{\|\mathrm Lu\|^2}-\,4\irdnu{\mathrm Lu:\mathrm M[u]}+\irdnu{\|\mathrm M[u]\|^2}\\
-\,8\irdnu{\(\mathrm Lu-\frac12\,\mathrm M[u]\):(\nabla u\otimes\nabla g)}\hspace*{2cm}\\
+\,2\irdnu{\left|\nabla u\right|^2\,(|\nabla g|^2-\Delta g)}-\,16\,\pi\,\lambda\irdeux{ \left|\nabla u \right|^2 e^{-u/2}}\,.
\end{multline*}
Since $\|\nabla u\otimes\nabla g\|^2=|\nabla u|^2\,|\nabla g|^2$, then we get for the corresponding trace free quantity
\[
\mathrm N u:=\nabla u\otimes\nabla g-\,\frac12\,(\nabla u\cdot\nabla g)\,\mathrm I_2
\]
that
\[
\|\,\mathrm Nu\|^2=|\nabla u|^2\,|\nabla g|^2-\,\frac12\,(\nabla u\cdot\nabla g)^2
\]
and hence obtain the identity
\begin{multline}\label{Eqno:EuclideanCondition}
0=4\irdnu{\left\|\,\mathrm Lu-\,\frac12\,\mathrm M[u]-\,\mathrm Nu\,\right\|^2}\\-\,2\irdeux{\left[\(\Delta g+ |\nabla g|^2-(\nabla g\cdot\omega)^2\)\,e^{-g} +\,8\,\pi\,\lambda\right]\,|\nabla u|^2 \, e^{-u/2}}
\end{multline}
where $\omega:=\nabla u/|\nabla u|$. If we assume that $\irdeux{|\nabla u|^2 \, e^{-u/2}}\neq0$ and define
\[
\Lambda:=\frac 1{8\,\pi}\,\frac{ \irdeux{\left[ (\nabla u\cdot\nabla g)^2-(\Delta g+|\nabla g|^2)\,|\nabla u|^2\right]
\, e^{-u/2-g}\,
}}{\irdeux{|\nabla u|^2 \, e^{-u/2}}}\,,
\]
then we get a contradiction if $\lambda<\Lambda$. To keep maximal generality, one could even include the term $\irdnu{\left\|\,\mathrm Lu-\,\frac12\,\mathrm M[u]-\,\mathrm Nu\,\right\|^2}$ in the definition of $\Lambda$, as it was done in the case of manifolds. However, this is a quite complicated criterion to verify since it involves the solution to~\eqref{eq:EL_weigthed} itself. Hence it makes sense to consider the simpler case where $\mu$ has radial symmetry. In that case it is also known from~\cite{MR634248} that $u$ is radially symmetric if $\mu$ is a monotone non-increasing function of $|x|$. Let
\[
\Lambda_\star:= - \frac 1{8\,\pi}\,\inf_{x\in\R^2}\(e^{-g}\,\Delta g\)=\inf_{x\in\R^2}\frac{-\,\Delta\log\mu}{8\,\pi\,\mu}\,.
\]
\begin{Thm}\label{Thm:Euclidean} Assume that $\mu$ is a radially symmetric function. Then any radially symmetric solution to~\eqref{eq:EL_weigthed} is a constant if $\lambda<\Lambda_\star$ and the inequality~\eqref{eq:weighted_onofri} holds with $\lambda=\Lambda_\star$ if equality is achieved among radial functions.\end{Thm}

\begin{Ex}\label{Ex1} The Euclidean Onofri inequality corresponds to
\[
\mu(x)=\frac1{\pi\,(1+|x|^2)^2}\quad\forall\,x\in\R^2\,.
\]
Since $-\Delta\log\mu=8\,\pi\,\mu$, it is known that $\Lambda_\star=1$ is the optimal constant. See~\cite{MR1143664} for further details. Let us notice that the analysis of the equation
\[
\mathrm Lu-\,\frac12\,\mathrm M[u]-\,\mathrm Nu=0
\]
in the case $\Lambda_\star=1$ provides a proof of the uniqueness of the radial solution to~\eqref{eq:EL_weigthed}, which is alternative to the result of~\cite{BKLN}.\end{Ex}

\begin{Ex}\label{Ex2} It is straightforward to deduce a perturbation result from Theorem~\ref{Thm:Euclidean}, that goes as follows. Let
\[
\mu(x)=\frac{e^{-h(x)}}{Z\,(1+|x|^2)^2}\quad\forall\,x\in\R^2\,,
\]
where $h$ is a radial function and $Z$ a normalization constant so that $\mu$ is a probability measure. We shall assume that~$h$ has a bounded variation and is such that $|x|^4\,\Delta h$ is bounded from below. Then we have the estimate
\[
\inf_{x\in\R^2}\frac{-\,\Delta\log\mu}{8\,\pi\,\mu}\ge\,e^{-\mathrm{Var}(h)}\,\left[1+\frac18\,\inf_{x\in\R^2}(1+|x|^2)^2\,\Delta h\right]\,.
\]
\end{Ex}

\begin{Ex}\label{Ex3} The subcritical Onofri inequality has been studied in~\cite{MR2996772}. It plays an important role for the study of the subcritical Keller-Segel model and its asymptotics for large times, and goes as follows. Let $\mu=n/M$ where $n$ is given as the unique (up to constants) radial solution to
\[
-\,\Delta c = n = M\,\frac{e^{c-\frac12|x|^2}}{\irdeux{e^{c-\frac12|x|^2}}}
\]
and the mass $M$ is taken in the interval $M\in(0,8\,\pi)$. According to~\cite{MR634248}, $n$ is radially symmetric as a consequence of the moving plane method. It is straightforward to check that
\[
\inf_{x\in\R^2}\frac{-\,\Delta\log\mu}{8\,\pi\,\mu}=\frac M{8\,\pi}+\inf_{x\in\R^2}\frac1{4\,\pi\,\mu}
\]
but the symmetry of the solution to~\eqref{eq:weighted_onofri} is not true for $\lambda>\frac M{8\,\pi}$ and $\frac M{8\,\pi}$ turns out to be the value of the optimal constant in~\eqref{eq:weighted_onofri}. See~\cite{Campos-CPDE,MR2996772} for further details.
\end{Ex}

\begin{Ex}\label{Ex4} The parabolic-parabolic Keller-Segel model has global in time solutions with mass larger than $8\,\pi$ for some values of its parameters, according to~\cite{springerlink:10.1007/s00285-010-0357-5}. The stationary solutions in self-similar variables can be written as
\[
-\,\Delta c = \varepsilon\,x\cdot\nabla c+n\quad\mbox{with}\quad n = M\,\frac{e^{c-\frac12|x|^2}}{\irdeux{e^{c-\frac12|x|^2}}}
\]
(where $\varepsilon>0$ is a given parameter) and have been shown to be radially symmetric in~\cite{NSY02}. To prove that a weighted Onofri inequality holds with $\mu=n/M$, it is therefore sufficient to establish the range of $\lambda\in(0,\Lambda_\star)$ such that the minimizer of $u\mapsto\irdeux{|\nabla u|^2}-\,16\,\pi\,\lambda\left[\log\(\irdmu{e^u}\)-\irdmu u\right]$ is radially symmetric, where
\[
\inf_{x\in\R^2}\frac{-\,\Delta\log\mu}{8\,\pi\,\mu}=\frac M{8\,\pi}+\inf_{x\in\R^2} \frac{\varepsilon\,x\cdot\nabla c +2}{8\,\pi\,\mu}\,.
\]
Such symmetry breaking issues are however known to be difficult: see for ins\-tance~\cite{0951-7715-27-3-435} for a discussion of a related problem.\end{Ex}

\section*{Acknowledgments} J.D.~and G.J.~have been supported by the projects \emph{STAB} and \emph{Kibord} of the French National Research Agency (ANR). J.D.~and M.J.E..~have been supported by the project \emph{NoNAP} of the French National Research Agency (ANR). The authors~thank M.~Loss and M.~Kowalczyk for many stimulating discussions.\\
{\sl\small\copyright~2016 by the authors. This paper may be reproduced, in its entirety, for non-commercial purposes.}

\end{document}